\documentclass[12pt,a4paper]{amsart}
\usepackage{amssymb}
\usepackage{mathrsfs}

\headsep=1cm \numberwithin{equation}{section}
\newtheorem{thm}{Theorem}[section]
\newtheorem{cor}[thm]{Corollary}
\newtheorem{lem}[thm]{Lemma}
\newtheorem{prop}[thm]{Proposition}
\theoremstyle{definition}
\newtheorem{defn}[thm]{Definition}
\theoremstyle{remark}
\newtheorem{rem}[thm]{Remark}
\numberwithin{equation}{section}

\begin{document}
\title{A note on the quintasymptotic prime ideals }
\author[Saeed Jahandoust and Reza Naghipour]{ Saeed Jahandoust and Reza Naghipour$^*$\\\\\\\,
\vspace*{0.5cm}Dedicated to Professor Peter Schenzel}
\address{Department of Mathematics, University of Tabriz, Tabriz, Iran.}
\email{saeed.e.jahan@gmail.com}
\address{Department of Mathematics, University of Tabriz, Tabriz, Iran;
and School of Mathematics, Institute for Research in Fundamental
Sciences (IPM), P.O. Box: 19395-5746, Tehran, Iran.}
\email{naghipour@ipm.ir} \email {naghipour@tabrizu.ac.ir}
\thanks{ 
This research was in part supported by a grant from IPM.\\
$^*$Corresponding author: e-mail: {\it naghipour@ipm.ir} (Reza Naghipour)}%

\begin{abstract}
Let $R$ denote a commutative  Noetherian ring, $I$ an ideal of $R$,
and let $S$ be a multiplicatively closed subset of $R$. In
\cite{Ra1}, Ratliff  showed that the sequence of sets ${\rm Ass}_RR/\overline{I}\subseteq {\rm Ass} _RR/\overline{I^2}
\subseteq {\rm Ass}_R R/\overline{I^3}\subseteq \dots $ increases and eventually stabilizes to a set  denoted $\overline{A^\ast}(I)$.
In \cite{Mc2}, S. McAdam  gave an interesting description of $\overline{A^\ast}(I)$ by making use of $R[It,t^{-1}]$, the Rees ring of $I$.
In this paper, we  give a second description of $\overline{A^\ast}(I)$ by making use of the Rees valuation rings of $I$. We also reprove a result
concerning when $\overline{I^n}R_S\cap R=\overline{I^n}$ for all integers $n>0$.

\end{abstract}
\maketitle
\section{Introduction}
Throughout this paper, all rings considered will be commutative and
Noetherian and will have non-zero identity elements. Such a ring will
be denoted by $R$, and the terminology is, in general, the same as that
in  \cite{BH}, \cite{Mc1} and \cite{Mc3}. Let $I$ be an ideal of $R$.
We denote by $\mathscr{R}$ the {\it graded Rees ring} $R[u,It]
:=\oplus_{n\in \mathbb{Z}}I^nt^n$ of $R$ with respect to  $I$, where
$t$ is an indeterminate and $u= t^{-1}$. Also, the integral closure of
$I$ in $R$ is denoted by $\overline{I}$, so
$\overline{I}:=\{x\in R\mid x\,\, \text{satisfies an equation of the form}\,\, x^n+c_1x^{n-1}+\dots+c_n=0, \text{where}\,\, c_i\in \,\,I^i \,\, \text{for}\,\, i=1,\dots, n\}.$ Further, if $(R, \frak m)$ is local, then $R^\ast$ denotes the completion of $R$ with respect to $\frak m$-adic topology.
The interesting concepts of quintasymptotic and asymptotic primes of $I$
were introduced and studied by McAdam in \cite{Mc2}, for generalizations of Schenzel's results in \cite{Sc}. These new
ideas have nice properties, for example see \cite{Ma}, \cite{Mc2} and \cite{Mc3}.
A prime ideal $\frak p$ of $R$ is called a {\it quintasymptotic} (resp.  {\it asymptotic}) prime
of an ideal $I$ of $R$, precisely when $I\subseteq \frak p$ and there exists a minimal prime $z$ in
 ${R_{\frak p}^\ast}$ with $\frak p{R_{\frak p}^\ast}$ minimal over $I{R_{\frak p}^\ast}+z$
(resp. $\frak p=\frak q\cap R$ for some quintasymptotic prime $\frak
q$ of the ideal $u\mathscr{R}$).

The set of {\it quintasymptotic} (resp. {\it asymptotic}) primes
of an ideal $I$ is denoted by $\overline{Q^\ast}(I)$ (resp. $\overline{A^{\ast\ast}}(I)$), so

\begin{center}
$\overline{Q^\ast}(I)= \{\frak p\in {\rm Spec} \,R\mid
I\subseteq \frak p$ and there is a minimal prime $z$ in ${R_\frak
p^\ast}$ with $\frak p{R_\frak p^\ast}$ minimal over $I{R_\frak
p^\ast} +z$ $\}$,  and\\

$\overline{A^{\ast\ast}}(I)= \{\frak p\cap R \mid \frak p \in
\overline{Q^\ast}(u\mathscr{R})\}$.
\end{center}
In \cite{Ra1}
L.J. Ratliff, Jr., showed that if $I$ is an ideal of a commutative
Noetherian ring $R$ with ${\rm ht}\, I>0,$ then the sequence of
associated prime ideals

$${\rm Ass}_R\, R/\overline{I^n},  n= 1,2, \ldots ,$$ becomes  eventually constant;
the stable value being denoted by $\overline{A^\ast}(I)$. Recall that later he has proved that the assumption
${\rm ht}\, I>0$ is not necessary.
As an interesting result, S. McAdam  in  \cite[Lemma 0.1]{Mc2}, showed that
$\overline{A^\ast}(I)=\overline{A^{\ast\ast}}(I).$ Using this result he
 proved that if $S$ is a multiplicatively closed subset of $R$ such that $S\subseteq
R\setminus\bigcup \{\frak q \cap R \mid \frak q\in
\overline{Q^\ast}(u\mathscr{R})\}$,  then $\overline{I^n}R_S\cap
R=\overline{I^n}$, for all $n>0$. This is a generalization of Schenzel's result in \cite[Theorem 5.6]{Sc}. In this paper we define the set
$\overline{B^\ast}(I)=\{\frak p\in {\rm Spec} \,R\mid
I\subseteq \frak p$ and there exists a minimal prime $z$ in $R$ with $z\subseteq \frak p$ and $\frak p/z$ is the center of a Rees
valuation ring of $I(R/z)$$\}$, and we will show that  $\overline{B^\ast}(I)=\overline{A^{\ast}}(I)$, by showing  $\overline{B^\ast}(I)=\overline{A^{\ast\ast}}(I).$ Then we will also use
$\overline{B^\ast}(I)=\overline{A^{\ast\ast}}(I)$ to reprove a result concerning those $Sُ$ for which  $\overline{I^n}R_S\cap
R=\overline{I^n}$, for all $n>0$.

The paper is organized as follows. In Section 2, we will prove  $\overline{B^\ast}(I)=\overline{A^\ast}(I)=\overline{A^{\ast\ast}}(I),$
and  in Section 3, which is the core of the paper, we will reprove \cite[Corollary 1.6]{Mc2}.

\section{$\overline{B^\ast}(I)=\overline{A^{\ast}}(I)$}
The purpose of this section is to establish $\overline{B^\ast}(I)=\overline{A^{\ast}}(I)$,
 by showing  $\overline{B^\ast}(I)=\overline{A^{\ast\ast}}(I).$
 The main goal of this section is Proposition 2.5,  which plays a key role in the proof of the main  theorem in Section 3.
 Let us, firstly,  recall some important results  concerning quintasymptotic primes and Rees valuation rings of an
ideal in a Noetherian ring.

\begin{lem}
Let $I$ be an ideal in a Noetherian ring $R$. Then:

{\rm(i) (see \cite[Lemma 2.1]{Mc2})} $\overline{Q^\ast}(I)$ is well defined finite set that contains all minimal prime divisors of
$I$.

{\rm (ii) (see \cite[Lemma 3.4]{Mc2})} $\frak p\in \overline{Q^\ast}(I)$ if and only if there is a
minimal prime $z$ of $R$ with $z\subseteq \frak p$ and $\frak p/z\in \overline{Q^\ast}(I(R/z))$.

{\rm (iii) ( see \cite[Proposition 3.8]{Mc2})} If $T$ is a finite module extension of $R$ such that $z\cap
T$ is a minimal prime in $R$ for all minimal primes  $z$ in $T$, then
$\overline{Q^\ast}(I)=\{\frak p\cap R \mid \frak p \in \overline{Q^\ast}(IT)\}$.
\end{lem}

\begin{lem}
Let $I$ be an ideal in a Noetherian ring $R$. Then $\frak p\in
\overline{A^{\ast\ast}}(I)$ if and only if there is a minimal prime
$z$ of $R$ with $z\subseteq \frak p$ and $\frak p/z \in
\overline{A^{\ast\ast}}(I(R/z))$.

\proof Let $\mathscr{R}$ be the Rees ring of $R$ with respect to $I$ and
$S=R[u,t]$. Then ${\rm Min}\,\mathscr{R}=\{zS\cap \mathscr{R} \mid
z\in{\rm Min}\,R\}$ by \cite[Theorem 1.5]{R}, and if $z$ is a
minimal prime of $R$, then $\mathscr{R}/{zS\cap \mathscr{R}}\cong
\mathscr{R}(R/z,I(R/z))$, by \cite[Lemma 1.1]{R}. Now, the
result follows from Lemma 2.1(ii). \qed \\
\end{lem}

\begin{defn}
Let $R$ be a Noetherian integral domain with the field of fractions
$K$ and let $I$ be an ideal of $R$. Let
$\mathscr{R}$ be the Rees ring of $R$ with respect
to $I$ and let $\mathscr{R}^\prime$ denote the integral closure of
$\mathscr{R}$ in the field of fractions of $\mathscr{R}$. For each
$x\in R$, let $V_I(x)$ be the largest positive integer such that
$x\in I^n$ (as usual, $I^0=R$ and $V_I(x)=\infty$ in case $x\in I^n$
for all $k\in \mathbb{N}$), and define,
$$\overline{V}_I(x)=\lim_{n\rightarrow\infty}\frac{V_I(x^n)}{n}.$$

It is known that:\\

{\rm(i)} $\overline{V}_I(x)$ is well defined for all $x\in
R$,  by \cite[Proposition 11.1]{Mc1}.

{\rm(ii)} Let $\frak{p}_1,\dots,\frak {p}_r$ be the
height one primes of $\mathscr{R}^\prime$ which contain $u$, and let
$v_i$ be the valuation associated with D.V.R.
$\mathscr{R}^\prime_{\frak{p}_i}$ and $e_i=v_i(u)$ for $i=1,\dots,r$.
Then for all $x\in R$, $\overline{V}_I(x)=\min\{\frac{v_i(x)}{e_i} |
i=1,2,\dots, r\}$, by \cite[Proposition 11.5]{Mc1}.

{\rm(iii)} For all positive integers $k$ and $x\in R$,
$\overline{V}_I(x)\geq k$ if and only if $x\in \overline{I^k}$ (as
usual, $\overline{I^0}=R$), by \cite[Corollary 11.6]{Mc1}.\\

The set of Rees valuation rings of $I$ is the set of rings
$V_i=\mathscr{R}^\prime_{\frak p_i} \cap K$, for $i=1,2, \dots, r$
and if $(V,\frak m_V)$ is a Rees valuation ring of $I$, then the center of
$V$ in $R$ is $\frak m_V\cap R$.
\end{defn}

\begin{rem}
Let $I$ be an ideal of a Noetherian ring $R$. Then:

{\rm(i)} An element $r\in R$ is in the integral closure of $I$ if
and only if, for every minimal prime $z$ in $R$, the image of $r$ in
$R/z$ is in the integral closure of $I(R/z)$.

{\rm(ii)} Let $z$ be a minimal prime ideal, and $\frak p$ a prime
ideal minimal over $I+z$. Then for all sufficiently large integers
$n$, for any ideal $J$ such that $I^n\subseteq J\subseteq
\overline{I^n}$, $\frak p$ is associated to $J$.

{\rm(iii)} With notations as in the Definition 2.3, if $R$ is integral domain, then
for all integers $n>0$,
$\overline{I^n}=\bigcap _{i=1}^r I^nV_i\cap R$ .
\end{rem}
\proof {\rm(i)} and {\rm(ii)} hold by \cite[Proposition 1.1.5]{HS} and
\cite[Lemma 5.4.4]{HS}. For prove {\rm(iii)}, if $I=0$ the result is clear. So
we assume that $I\neq 0$. Then, since $I^nV_i$ is principal and
$V_i$ is integrally closed, it follows that $\overline{I^n}\subseteq
\overline{I^n}V_i\subseteq \overline{I^nV_i}=I^nV_i$ for $i=1,2,
\dots, r$, and so $\overline{I^n}\subseteq\bigcap _{i=1}^r
I^nV_i\cap R$. Now, let $x\in I^nV_i\cap R$ for all $1\leq i\leq r$.
By \cite[Proposition 3.6]{Ra}, $It\subseteq \mathscr{R}^\prime
\setminus \frak p_i$,  so $It\subseteq \mathscr{R}^\prime _{\frak
p_i}\setminus {\frak p_i}\mathscr{R}^\prime _{\frak p_i}$. Thus
$v_i(u)=v_i(I)(:=\min\{v_i(a) \mid a\in I \})$.  Therefore
$v_i(x)\geq nv_i(u)$, and  hence $x\in \overline{I^n}$, by the Definition 2.3. \qed\\

Now we are prepared to prove the main result of this section, which plays a key role in the proof of the main  theorem in Section 3.

\begin{prop}
Let $I$ be an ideal in a Noetherian ring $R$. Then $\overline{B^\ast}(I)=\overline{A^{\ast\ast}}(I)$.
\end{prop}

\proof By Lemma 2.2, we may assume that $R$ is an integral domain with the field
of fraction $K$. Let $\frak p \in \overline{A^{\ast \ast}}(I)$ and suppose that
$\mathscr{R}$ is the Rees ring of $R$ with respect to $I$ with
the integral closure $\mathscr{R}^\prime$ in its field of fractions.
There exists $\frak q\in \overline{Q^\ast}(u\mathscr{R})$ such that
$\frak p=\frak q\cap R$. Put $\mathscr{S}=\mathscr{R}_{\frak q}$.
Then there exists a minimal prime $w$ in $\mathscr{S}^\ast$ such
that $u\mathscr{S}^\ast+w$ is a $\frak q\mathscr{S}^\ast$-primary
ideal. Since, $u^n\mathscr{S}^\ast \subseteq
\overline{u^n\mathscr{S}}\mathscr{S}^\ast\subseteq
\overline{u^n\mathscr{S}^\ast}$, it follows from the Remark 2.4(ii) that, for all sufficiently
large $n$, $\frak q\mathscr{S}^\ast$ is an associated prime to
$\overline{u^n\mathscr{S}}\mathscr{S}^\ast$ and so $\frak
q\mathscr{S}$ is an associated prime to $\overline{u^n\mathscr{S}}$. Thus
$\frak q$ is an associated prime to $\overline{u^n\mathscr{R}}$. Moreover, as
$\overline{u^n\mathscr{R}}=u^n\mathscr{R}^\prime\cap \mathscr{R}$,
there exists a prime ideal $\frak P$ in $\mathscr{R}^\prime$
such that $\frak q=\frak P\cap \mathscr{R}$ and $\frak P$ is
associated prime to $u^n\mathscr{R}^\prime$. Therefore, it is easy to
check that $\frak p$ is the center of Rees valuation ring $\frak
P\mathscr{R}^\prime_\frak P\cap K$ of $I$. That is $\frak p\in \overline{B^\ast}(I)$, and
so $\overline{A^{\ast\ast}}(I)\subseteq\overline{B^\ast}(I)$.

Now, let $\frak p\in \overline{B^\ast}(I)$. Then
there exists a prime divisor $\frak P$ of $u\mathscr{R}^\prime$,
such that $\frak p=\frak P\cap R$. Moreover, by \cite[Lemma
4.8.4]{HS}, there exists a ring
$\mathscr{R}\subseteq\mathscr{T}\subseteq\mathscr{R}^\prime$ and a
prime ideal $\frak q$ in $\mathscr{T}$ such that $\frak P$ contracts
to $\frak q$ and ${\rm ht}\frak P={\rm ht}\frak q$.  Also, by the proof
of \cite[Lemma 4.8.4]{HS}, $\mathscr{T}$ is a finite
$\mathscr{R}$-module. Since $u\mathscr{T}\subseteq \frak q$ and ${\rm ht}\frak q=1$,
it follows that $\frak q\in
\overline{Q^{\ast}}(u\mathscr{T})$ by Lemma 2.1(i).  Hence, by Lemma 2.1(iii), we have
$$\frak q\cap \mathscr{R}=(\frak P\cap \mathscr{T})\cap
\mathscr{R}=\frak P\cap \mathscr{R} \in
\overline{Q^{\ast}}(u\mathscr{R}).$$ Thus $(\frak P\cap
\mathscr{R})\cap R=\frak p\in \overline{A^{\ast\ast}}(I)$, and so $\overline{B^\ast}(I)\subseteq \overline{A^{\ast\ast}}(I)$.
This completes the proof. \qed \\

\begin{cor}
Let $I$ be an ideal in a Noetherian ring $R$. Then $\overline{B^\ast}(I)=\overline{A^{\ast}}(I)$.
\end{cor}

\proof The assertion follows from Proposition 2.5 and \cite[Lemma 0.1]{Mc2}. \qed \\

\section{McAdam's question}
McAdam, in \cite[Corollary 1.6]{Mc2} pointed out the following trivial consequence of his result that
$\overline{A^\ast}(I)=\overline{A^{\ast\ast}}(I)$.  If $S$ is a
multiplicatively closed subset of $R\setminus\bigcup \{\frak p \mid \frak p\in \overline{A^{\ast\ast}}(I)\}$,
then  for all  $n\geq1$,  $\overline{I^n}R_S\cap R=\overline{I^n}$. He asked if that could be proven without
using that $\overline{A^\ast}(I)=\overline{A^{\ast\ast}}(I)$. In this section we show that it can be, using
$\overline{A^{\ast\ast}}(I)=\overline{B^\ast}(I)$. 

\begin{thm}
Let $I$ be an ideal in a Noetherian ring $R$ and let $S$ be a
multiplicatively closed subset of $R$ such that $S\subseteq R\setminus
\bigcup \{\frak p \mid \frak p\in \overline{A^{\ast\ast}}(I)\}$.
Then $\overline{I^n}R_S\cap R=\overline{I^n}$ for all integers $n>0$.
\end{thm}

\proof By Lemma 2.2 and Remark 2.4(i), we may assume that $R$ is an integral domain.
Let $n$ be a positive integer and let $x\in \overline{I^n}R_S\cap R$.
Let $(V,\frak m_V)$ be a Rees valuation ring of $I$. Then there
exists $s\in S$ such that $xs\in \overline{I^n}$, and hence $xs\in
I^nV$ by the Remark 2.4(iii). Now, in view of the Proposition 2.5, $s$ is not in any center of Rees valuation
rings of $I$.  In particular, $s$ is not in $\frak m_V$, and so $x\in
I^nV$. Thus $x\in \overline{I^n}$, by the Remark 2.4(iii). As the opposite inclusion is obvious,
the result follows. \qed \\

\begin{center}
{\bf Acknowledgments}
\end{center}
The authors are deeply grateful to the referee for  careful reading and many useful suggestions.
The authors  would  like to thank  Professor L.J. Ratliff for
his careful reading of the first draft and kind comments in the preparation of this article.  Also, we would
like to thank from the School of Mathematics, Institute for Research in Fundamental
Sciences (IPM) for its financial support.


\end{document}